# BREAKDOWN POINTS FOR MAXIMUM LIKELIHOOD ESTIMATORS OF LOCATION–SCALE MIXTURES

BY CHRISTIAN HENNIG

*Universität Hamburg*

ML-estimation based on mixtures of Normal distributions is a widely used tool for cluster analysis. However, a single outlier can make the parameter estimation of at least one of the mixture components break down. Among others, the estimation of mixtures of $t$-distributions by McLachlan and Peel [*Finite Mixture Models* (2000) Wiley, New York] and the addition of a further mixture component accounting for "noise" by Fraley and Raftery [*The Computer J.* **41** (1998) 578–588] were suggested as more robust alternatives. In this paper, the definition of an adequate robustness measure for cluster analysis is discussed and bounds for the breakdown points of the mentioned methods are given. It turns out that the two alternatives, while adding stability in the presence of outliers of moderate size, do not possess a substantially better breakdown behavior than estimation based on Normal mixtures. If the number of clusters $s$ is treated as fixed, $r$ additional points suffice for all three methods to let the parameters of $r$ clusters explode. Only in the case of $r = s$ is this not possible for $t$-mixtures. The ability to estimate the number of mixture components, for example, by use of the Bayesian information criterion of Schwarz [*Ann. Statist.* **6** (1978) 461–464], and to isolate gross outliers as clusters of one point, is crucial for an improved breakdown behavior of all three techniques. Furthermore, a mixture of Normals with an improper uniform distribution is proposed to achieve more robustness in the case of a fixed number of components.

**1. Introduction.** ML-estimation based on mixtures of Normal distributions (NMML) is a flexible and widely used technique for cluster analysis [e.g., Wolfe (1967), Day (1969), McLachlan (1982), McLachlan and Basford (1988), Fraley and Raftery (1998) and Wang and Zhang (2002)]. Moreover, it is applied to density estimation and discrimination [Hastie and Tibshirani (1996) and Roeder and Wasserman (1997)]. Banfield and Raftery (1993) introduced the term "model-based cluster analysis" for such methods.

Observations $x_1, \ldots, x_n$ are modeled as i.i.d. with density

$$f_\eta(x) = \sum_{j=1}^{s} \pi_j \varphi_{a_j, \sigma_j^2}(x), \tag{1.1}$$

where $\eta = (s, a_1, \ldots, a_s, \sigma_1, \ldots, \sigma_s, \pi_1, \ldots, \pi_s)$ is the parameter vector, the number of components $s \in \mathbb{N}$ may be known or unknown, $a_j \in \mathbb{R}$, $\sigma_j > 0$, $\pi_j \geq 0$, $j = 1, \ldots, s$, $\sum_{j=1}^{s} \pi_j = 1$ and $\varphi_{a,\sigma^2}$ denotes the density of a Normal distribution with mean $a$ and variance $\sigma^2$, $\varphi = \varphi_{0,1}$. Mixtures of multivariate Normals are often used, but for the sake of simplicity, considerations are restricted to the case of one-dimensional data in this paper. The results essentially carry over to the multivariate case.

As in many other ML-techniques that are based on the Normal distribution, NMML is not robust against gross outliers, in particular, if the number of components $s$ is treated as fixed: the estimators of the parameters $a_1, \ldots, a_s$ are weighted means of the observations. For each observation, the weights sum up to 1 [see Redner and Walker (1984)], which means that at least one of these parameters can become arbitrarily large if a single extreme point is added to a dataset.

There are some ideas in the literature to overcome the robustness problems of Normal mixtures. The software MCLUST [Fraley and Raftery (1998)] allows the addition of a mixture component accounting for "noise," modeled as a uniform distribution on the convex hull (the range in one dimension, respectively) of the data. The software EMMIX [Peel and McLachlan (2000)] can be used to fit a mixture of $t$-distributions instead of Normals. Further, it has been proposed to estimate the component parameters by more robust estimators [Campbell (1984), McLachlan and Basford (1988) and Kharin (1996), page 275], in particular, by Huber's (1964, 1981) M-estimators corresponding to ML-estimation for a mixture of Huber's least favorable distributions [Huber (1964)].

While a clear gain of stability can be demonstrated for these methods in various examples [see, e.g., Banfield and Raftery (1993) and McLachlan and Peel (2000), page 231 ff.], there is a lack of theoretical justification of their robustness. Only Kharin [(1996), page 272 ff.] obtained some results for fixed $s$. He showed that under certain assumptions on the speed of convergence of the proportion of contamination to 0 with $n \to \infty$, Huber's M-estimation is asymptotically superior to NMML. In the present paper, mixtures of a class of location–scale models are treated including the aforementioned distributions. The addition of a "noise" component is also investigated.

Up to now there is no agreement about adequate robustness measures for cluster analysis. In a model-based cluster analysis, the clusters are characterized by the parameters of their mixture components. For fixed $s$ an influence function [Hampel (1974)] and a breakdown point [Hampel (1971) and Donoho and Huber (1983)] for these parameters can easily be defined. The "addition breakdown point" is the minimal proportion of points to be added to an original dataset so that the parameter estimator for the new dataset deviates as far as possible from the value obtained from the original dataset. However, there are some particular issues in cluster analysis. Partitioning methods may possess a bounded influence function and the minimal possible breakdown point at the same time. The breakdown point may depend strongly on the constellation of the data points [Garcia-Escudero and Gordaliza (1999)]. One may distinguish between breakdown of a single cluster and breakdown of all clusters [Gallegos (2003)], and breakdown could be characterized by means of the classification of the points instead of the estimated parameters [Kharin (1996), page 49]. The breakdown concepts in the literature cited above only apply to a fixed number of components $s$. If $s$ is estimated, there are data constellations "on the border" between two different numbers of components, leading to different numbers of parameters to estimate.

The outline of the paper is as follows. In Section 2, the techniques treated in this paper and their underlying models are introduced.

In Section 3, robustness measures and breakdown points in terms of parameters (Definition 3.1 for fixed $s$, Definition 3.2 for estimated $s$) as well as of classification (Definition 3.4) are defined.

In Section 4, results about the parameter breakdown behavior of the mixture-based clustering techniques are derived. It is shown that all discussed techniques have a breakdown point of $r/(n+r)$ for $r < s$ of the mixture components in the case of fixed $s$ (Theorem 4.4). A better breakdown behavior can be attained by maximizing a kind of "improper likelihood" where "noise" is modeled by an improper uniform distribution on the real line (Theorem 4.11). For the case of estimated $s$, using an information criterion [Akaike (1974) and Schwarz (1978)], a breakdown point larger than $1/(n+1)$ can be attained for all considered methods. They all are able to isolate gross outliers as new mixture components on their own and are therefore very stable against extreme outliers. However, breakdown can happen because additional points inside the area of the estimated mixture components of the original data can lead to the estimation of a smaller number of components (Theorems 4.13 and 4.16). Some numerical examples are given, illustrating the relative stability of the methods and the nonequivalence of parameter and classification breakdown and of addition and replacement breakdown. Some data constellations turn out to be so stable that they lead to an addition parameter breakdown point larger than $1/2$. The paper is completed by some concluding discussions.

**2. Models and methods.** The Normal mixture (1.1) belongs to the class of mixtures of location–scale families $f_\eta$ which can be defined as follows:

$$(2.1) \qquad f_\eta(x) = \sum_{j=1}^{s} \pi_j f_{a_j,\sigma_j}(x), \qquad \text{where } f_{a,\sigma}(x) = \frac{1}{\sigma} f\left(\frac{x-a}{\sigma}\right),$$

where $\eta$ is defined as in (1.1). Assume that

$(2.2) \qquad\qquad\qquad f$ is symmetrical about $0$,

$(2.3) \qquad\qquad\qquad f$ decreases monotonicly on $[0, \infty]$,

$(2.4) \qquad\qquad\qquad f > 0$ on $\mathbb{R}$,

$(2.5) \qquad\qquad\qquad f$ is continuous.

Besides the $\mathcal{N}(0,1)$-distribution, these assumptions are fulfilled, for example, for the $t_\nu$-distribution with $\nu$ degrees of freedom and for Huber's least favorable distribution, used as a basis for mixture modeling in Peel and McLachlan (2000) and McLachlan and Basford (1988), respectively.

The following properties will be referred to later. It follows from (2.2)–(2.4) that, for given points $x_1, \ldots, x_n$ and a compact set $C = [a,b] \times [\mu, \xi] \subset \mathbb{R} \times \mathbb{R}^+$ (this notation implies $\mu > 0$ here),

$$(2.6) \qquad \inf\{f_{a,\sigma}(x) : x \in \{x_1, \ldots, x_n\}, (a,\sigma) \in C\} = f_{\min} > 0.$$

For fixed $x$, $\lim_{m\to\infty} a_m = \infty$ and arbitrary sequences $(\sigma_m)_{m \in \mathbb{N}}$, observe that

$$(2.7) \qquad \lim_{m\to\infty} f_{a_m,\sigma_m}(x) \leq \lim_{m\to\infty} \min\left(\frac{1}{\sigma_m} f(0), \frac{1}{\sigma_0} f\left(\frac{x-a_m}{\sigma_m}\right)\right) = 0$$

as long as $\sigma_m \geq \sigma_0 > 0$.

The addition of a uniform mixture component on the range of the data is also considered, which is the one-dimensional case of a suggestion by Banfield and Raftery (1993). That is, for given $x_{\min} < x_{\max} \in \mathbb{R}$,

$$(2.8) \qquad f_\zeta(x) = \sum_{j=1}^{s} \pi_j f_{a_j,\sigma_j}(x) + \pi_0 \frac{\mathbb{1}(x \in [x_{\min}, x_{\max}])}{x_{\max} - x_{\min}},$$

where $\zeta = (s, a_1, \ldots, a_s, \sigma_1, \ldots, \sigma_s, \pi_0, \pi_1, \ldots, \pi_s)$, $\pi_0, \ldots, \pi_s \geq 0$, $\sum_{j=0}^{s} \pi_j = 1$ and $\mathbb{1}(A)$ is the indicator function for the statement $A$.

The log-likelihood functions for the models (2.1) and (2.8) for given data $\mathbf{x}_n$, with minimum $x_{\min,n}$ and maximum $x_{\max,n}$ (this notation is also used later), are

$$(2.9) \qquad L_{n,s}(\eta, \mathbf{x}_n) = \sum_{i=1}^{n} \log\left(\sum_{j=1}^{s} \pi_j f_{a_j,\sigma_j}(x_i)\right),$$

$$(2.10) \qquad L_{n,s}(\zeta, \mathbf{x}_n) = \sum_{i=1}^{n} \log\left(\sum_{j=1}^{s} \pi_j f_{a_j,\sigma_j}(x_i) + \frac{\pi_0}{x_{\max,n} - x_{\min,n}}\right).$$

As can easily be seen by setting $a_1 = x_1$, $\sigma_1 \to 0$, $L_{n,s}$ can become arbitrarily large for $s > 1$. Thus, to define a proper ML-estimator, the parameter space must be suitably restricted. The easiest restriction is to specify $\sigma_0 > 0$ and to demand

$$\sigma_j \geq \sigma_0, \qquad j = 1, \ldots, s. \tag{2.11}$$

This is used, for example, in DeSarbo and Cron (1988) and may easily be implemented with the EM-algorithm [Dempster, Laird and Rubin (1977) and Redner and Walker (1984); see Lemma 2.1], the most popular routine to compute mixture ML-estimators. A drawback of this restriction is that the resulting ML-estimators are no longer scale equivariant because the scale of the data can be made smaller than $\sigma_0$ by multiplication with a constant. The alternative restriction

$$\min_{j,k=1,\ldots,s} \sigma_j/\sigma_k \geq c \tag{2.12}$$

for fixed $c \in (0, 1]$ leads to properly defined, scale-equivariant, consistent ML-estimators for the Normal case $f = \varphi_{0,1}$ without noise [Hathaway (1985)]. This includes the popular simplification $\sigma_1 = \cdots = \sigma_s$, which corresponds to $k$-means clustering and is the one-dimensional case of some of the covariance parameterizations implemented in MCLUST [Fraley and Raftery (1998)]. However, unless $c = 1$, the computation is not straightforward [Hathaway (1986)]. Furthermore, the restriction (2.12) cannot be applied to the model (2.8), because the log-likelihood function may be unbounded; see Lemma A.1. For the case of fixed $s$, Corollary 4.5 says that estimation based on (2.12) does not yield better breakdown properties than its counterpart using (2.11). Therefore, the restriction (2.11) is used for all other results. Guidelines for the choice of $\sigma_0$ and $c$ are given in Section A.1. For results about consistency of local maximizers of the log-likelihood function, see Redner and Walker (1984).

The following lemma summarizes some useful properties of the maximum likelihood estimators, which follow from the derivations of Redner and Walker (1984).

NOTATION. Let $\theta_j = (a_j, \sigma_j)$, $j = 1, \ldots, s$, $\theta = (\theta_1, \ldots, \theta_s)$ denote the location and scale parameters of $\eta$ and $\zeta$, respectively, $\theta^*$, $\eta^*$, $\zeta^*$ by analogy. The parameters included in $\eta^*$, $\zeta^*$ will be denoted by $s^*$, $a_1^*$, $\pi_1^*$ and so on, and by analogy for $\hat{\eta}$, $\hat{\zeta}, \ldots$.

LEMMA 2.1. *For given $\eta$, let*

$$p_{ij} = \frac{\pi_j f_{a_j,\sigma_j}(x_i)}{\sum_{k=1}^s \pi_k f_{a_k,\sigma_k}(x_i)}, \qquad i = 1, \ldots, n. \tag{2.13}$$

*A maximizer* $\hat{\eta}$ *of*

$$\text{(2.14)} \quad \sum_{j=1}^{s}\left[\sum_{i=1}^{n}p_{ij}\log\pi_j^*\right]+\sum_{j=1}^{s}\sum_{i=1}^{n}p_{ij}\log f_{a_j^*,\sigma_j^*}(x_i)$$

*over* $\eta^*$ *leads to an improvement of* $L_{n,s}$ *unless* $\eta$ *itself attains the maximum of* (2.14).

*For given* $\zeta$ *in* (2.10), *the same statements hold with*

$$\text{(2.15)} \quad \begin{aligned} p_{ij} &= \frac{\pi_j f_{a_j,\sigma_j}(x_i)}{\sum_{k=1}^{s}\pi_k f_{a_k,\sigma_k}(x_i)+\pi_0/(x_{\max,n}-x_{\min,n})}, \quad j=1,\ldots,s, \\ p_{i0} &= \frac{\pi_0/(x_{\max,n}-x_{\min,n})}{\sum_{k=1}^{s}\pi_k f_{a_k,\sigma_k}(x_i)+\pi_0/(x_{\max,n}-x_{\min,n})}. \end{aligned}$$

*In* (2.14), *the first sum starts at* $j=0$.

*For any global maximizer* $\eta$ *as well as* $\zeta$ *of* $L_{n,s}$ *for given* $\mathbf{x}_n$, *the following conditions hold under* (2.11) *for* $j=1,\ldots,s$ *with* $p_{ij}$, $i=1,\ldots,n$:

$$\text{(2.16)} \quad \pi_j = \frac{1}{n}\sum_{i=1}^{n}p_{ij},$$

$$\text{(2.17)} \quad \begin{aligned} (a_j,\sigma_j) &= \arg\max S_j(a_j^*,\sigma_j^*) \\ &= \arg\max \sum_{i=1}^{n}p_{ij}\log\left(\frac{1}{\sigma_j^*}f\left(\frac{x_i-a_j^*}{\sigma_j^*}\right)\right). \end{aligned}$$

*In case of* (2.10), *property* (2.16) *holds for* $j=0$ *as well.*

Note that (2.13) defines the so-called E-step, and maximization of (2.14) defines the so-called M-step of the EM-algorithm, where the two steps are alternately carried out.

LEMMA 2.2. *Under* (2.11), *with*

$$C = [x_{\min,n}, x_{\max,n}] \times \left[\sigma_0, \frac{\sigma_0 f(0)}{f((x_{\max,n}-x_{\min,n})/\sigma_0)}\right]$$

*and* $\pi_1,\ldots,\pi_s > 0$,

$$\text{(2.18)} \quad \forall\, \theta^* \notin C^s \; \exists\, \theta \in C^s : L_{n,s}(\eta) > L_{n,s}(\eta^*).$$

Proofs are given in Section A.2.

Note that $L_{n,s}$ is continuous [cf. (2.5)] and a global maximizer has to lie in $C^s \times [0,1]^s$ because of (2.18). Therefore, we have following result.

COROLLARY 2.3. *Under the restriction* (2.11), *there exists a* (*not necessarily unique*) *global maximum of* $L_{n,s}$ *with arguments in* $C^s \times [0,1]^s$.

For NMML and (2.12), this is shown by Hathaway (1985). Define $\eta_{n,s} = \arg\max L_{n,s}$ and $\zeta_{n,s}$ analogously. In the case of nonuniqueness, $\eta_{n,s}$ can be defined as an arbitrary maximizer, for example, the lexicographically smallest one. The $p_{ij}$-values from (2.13) and (2.15), respectively, can be interpreted as the a posteriori probabilities that a point $x_i$ had been generated by component $j$ under the a priori probability $\pi_j$ for component $j$ with parameters $a_j, \sigma_j$. These values can be used to classify the points and to generate a clustering by

$$(2.19) \qquad l(x_i) = \arg\max_j p_{ij}, \qquad i = 1, \ldots, n,$$

where the ML-estimator is plugged into the definition of $p_{ij}$.

All theorems derived in the present paper will hold for any of the maximizers. For ease of notation, $\eta_{n,s}$ and $\zeta_{n,s}$ will be treated as well defined in the following. Note that, for $s > 1$, nonuniqueness always occurs due to "label switching" of the mixture components. Further, for ease of notation, it is not assumed, in general, that $\pi_j > 0 \ \forall j$ or that all $(a_j, \sigma_j)$ are pairwise distinct.

Consider now the number of mixture components $s \in \mathbb{N}$ as unknown. The most popular method to estimate $s$ is the use of information-based criteria such as AIC [Akaike (1974)] and BIC [Schwarz (1978)]. The latter is implemented in MCLUST. EMMIX computes both. The estimator $s_n$ for the correct order of the model is defined as $s_n = \arg\max_s C(s)$, where

$$(2.20) \qquad \begin{aligned} C(s) &= \mathrm{AIC}(s) = 2L_{n,s}(\eta_{n,s}) - 2k \quad \text{or} \\ C(s) &= \mathrm{BIC}(s) = 2L_{n,s}(\eta_{n,s}) - k\log n, \end{aligned}$$

where $k$ denotes the number of free parameters, that is, $k = 3s - 1$ for (2.1) and $k = 3s$ for (2.8). Under assumptions satisfied under (2.11) but not under (2.12) for the models discussed here (compare Lemma A.1), Lindsay [(1995), page 22] shows that the number of distinct points in the dataset is an upper bound for the maximization of $L_{n,s}(\eta_{n,s})$ over $s$, and therefore for the maximization of $C(s)$ as well. Thus, only a finite number of values for $s$ have to be investigated to maximize $C(s)$ and this means that (again not necessarily unique) maximizers exist.

While the AIC is known to overestimate $s$ asymptotically [see, e.g., Bozdogan (1994)], the BIC is shown at least in some restricted situations to be consistent in the mixture setup [Keribin (2000)]. I mainly consider the BIC here. Further suggestions to estimate $s$, which are more difficult to analyze with respect to the breakdown properties, are given, for example, by Bozdogan (1994) and Celeux and Soromenho (1996). EMMIX also allows the estimation of $s$ via a bootstrapped likelihood ratio test [McLachlan (1987)].

**3. Breakdown measures for cluster analysis.** The classical meaning of "addition breakdown" for finite samples is that an estimator can be driven arbitrarily far away from its original value by addition of unfortunate data points, usually by gross outliers. For "replacement breakdown points", points from the original sample are replaced [Donoho and Huber (1983)]. Zhang and Li (1998) and Zuo (2001) derive relations between these two concepts. In the present paper, addition breakdown is considered. Breakdown means that estimators that can take values on the whole range of $\mathbb{R}^p$ can leave every compact set. If the range of values of a parameter is bounded, breakdown means that the addition of points can take the estimator arbitrarily close to the bound, for example, a scale parameter to 0. Such a definition is relatively easily applied to the estimation of mixture components, but it cannot be used to compare the robustness of mixture estimators with other methods of cluster analysis.

Therefore, the more familiar parameter breakdown point will be defined first. Then, a breakdown definition in terms of the classification of points to clusters is proposed.

A "parameter breakdown" can be understood in two ways. A situation where at least one of the mixture components explodes is defined as breakdown in Garcia-Escudero and Gordaliza (1999). That is, breakdown occurs if the whole parameter vector leaves all compact sets [not including scales of 0 under (2.12)]. In contrast, Gallegos (2003) defines breakdown in cluster analysis as a situation where *all* clusters explode simultaneously. Intermediate situations may be of interest in practice, especially if a researcher tries to prevent the breakdown of a single cluster by specifying the number of clusters to be larger than expected, so that additional clusters can catch the outliers. This is discussed (but not recommended—in agreement with the results given here) by Peel and McLachlan (2000). The definition given next is flexible enough to account for all mentioned situations.

DEFINITION 3.1. Let $(E_n)_{n \in \mathbb{N}}$ be a sequence of estimators of $\eta$ in model (2.1) or of $\zeta$ in model (2.8) on $\mathbb{R}^n$ for fixed $s \in \mathbb{N}$. Let $r \leq s$, $\mathbf{x}_n = (x_1, \ldots, x_n)$ be a dataset, where

(3.1) $$\hat{\eta} = E_n(\mathbf{x}_n) \Rightarrow \hat{\pi}_j, \qquad j = 1, \ldots, s.$$

The *r-components parameter breakdown point* of $E_n$ is defined as

$$B_{r,n}(E_n, \mathbf{x}_n) = \min_g \bigg\{ \frac{g}{n+g} : \exists j_1 < \cdots < j_r$$
$$\forall D = [\pi_{\min}, 1] \times C, \ \pi_{\min} > 0,$$
$$C \subset \mathbb{R} \times \mathbb{R}^+ \text{ compact } \exists \mathbf{x}_{n+g} = (x_1, \ldots, x_{n+g}),$$
$$\hat{\eta} = E_{n+g}(\mathbf{x}_{n+g}) : (\hat{\pi}_j, \hat{a}_j, \hat{\sigma}_j) \notin D, j = j_1, \ldots, j_r \bigg\}.$$

The proportions $\hat{\pi}_j$ are defined not to break down if they are bounded away from 0, which implies that they are bounded away from 1 if $s > 1$. Assumption (3.1) is necessary for the definition to make sense; $\hat{\pi}_j = 0$ would imply that the corresponding location and scale parameters could be chosen arbitrarily large without adding any point. Condition (3.1) may be violated for ML-estimators in situations where $s$ is not much smaller than $n$, but these situations are usually not of interest in cluster analysis. In particular, (3.1) does not hold if $s$ exceeds the number of distinct $x_i$; see Lindsay [(1995), page 23].

The situation of $\hat{\pi}_0 \to 0$ in model (2.8) is not defined as breakdown, because the noise component is not considered as an object of interest in itself in this setup.

In the case of unknown $s$, considerations are restricted to the case of one-component breakdown. Breakdown robustness means that neither of the $s$ mixture components estimated for $\mathbf{x}_n$ vanishes, nor that any of their scale and location parameters explodes to $\infty$ under addition of points. It is, however, allowed that the new dataset yields more than $s$ mixture components and that the additional mixture components have arbitrary parameters. This implies that, if the outliers form a cluster on their own, their component can simply be added without breakdown. Further, breakdown of the proportions $\pi_j$ to 0 is no longer of interest when estimating $s$ according to the AIC or BIC, because if some $\pi_j$ is small enough, component $j$ can be simply left out, and the other proportions can be updated to sum up to 1. This solution with $s - 1$ clusters leads approximately to the same log-likelihood and will be preferred due to the penalty on the number of components:

DEFINITION 3.2. Let $(E_n)_{n \in \mathbb{N}}$ be a sequence of estimators of $\eta$ in model (2.1) or of $\zeta$ in model (2.8) on $\mathbb{R}^n$, where $s \in \mathbb{N}$ is estimated as well. Let $\mathbf{x}_n = (x_1, \ldots, x_n)$ be a dataset. Let $s^*$ be the estimated number of components of $E_n(\mathbf{x}_n)$. The *parameter breakdown point* of $E_n$ is defined as

$$B_n(E_n, \mathbf{x}_n) = \min_g \left\{ \frac{g}{n+g} : \forall C \subset \mathbb{R}^s \times (R^+)^s \text{ compact} \right.$$
$$\exists \mathbf{x}_{n+g} = (x_1, \ldots, x_{n+g}), \ \hat{\eta} = E_{n+g}(\mathbf{x}_{n+g}):$$
$$\text{pairwise distinct } j_1, \ldots, j_{s^*} \text{ do not exist,}$$
$$\left. \text{such that } (\hat{a}_{j_1}, \ldots, \hat{a}_{j_{s^*}}, \hat{\sigma}_{j_1}, \ldots, \hat{\sigma}_{j_{s^*}}) \in C \right\}.$$

This implies especially that breakdown occurs whenever $\hat{s} < s^*$.

Now, the classification breakdown is defined. A mapping $E_n$ is called a general clustering method (GCM) if it maps a set of entities $\mathbf{x}_n = \{x_1, \ldots, x_n\}$ to a collection of subsets $\{C_1, \ldots, C_s\}$ of $\mathbf{x}_n$. A special case are partitioning methods where $C_i \cap C_j = \varnothing$ for $i \neq j \leq s$, $\bigcup_{j=1}^{s} C_s = \mathbf{x}_n$. An ML-mixture estimator induces a partition by (2.19) and $C_j = \{x_i : l(x_i) = j\}$, given a rule to break ties in the $p_{ij}$.

If $E_n$ is a GCM and $\mathbf{x}_{n+g}$ is generated by adding $g$ points to $\mathbf{x}_n$, $E_{n+g}(\mathbf{x}_{n+g})$ induces a clustering on $\mathbf{x}_n$, which is denoted by $E_n^*(\mathbf{x}_{n+g})$. Its clusters are denoted by $C_1^*, \ldots, C_{s^*}^*$. If $E_n$ is a partitioning method, $E_n^*(\mathbf{x}_{n+g})$ is a partition as well. Note that $s^*$ may be smaller than $s$ when $E_n$ produces $s$ clusters for all $n$. Assume in the following that $E_n$ is a partitioning method. The resulting definition may be tentatively applied to other clustering methods as well.

As will be illustrated in Remark 4.18, different clusters of the same data may have a different stability. Thus, it makes sense to define robustness with respect to the individual clusters. This requires a measure for the similarity between a cluster of $E_n^*(\mathbf{x}_{n+g})$ and a cluster of $E_n(\mathbf{x}_n)$, that is, between two subsets $C$ and $D$ of some finite set. The following proposal equals 0 only for disjoint sets and 1 only for equal sets:

$$\gamma(C, D) = \frac{2|C \cap D|}{|C| + |D|}.$$

The definition of (addition) breakdown is based on the similarity of a cluster $C \in E_n(\mathbf{x}_n)$ to its most similar cluster in $E_n^*(\mathbf{x}_{n+g})$. A similarity between $C$ and a partition $\hat{E}_n$ is defined by

$$\gamma^*(C, \hat{E}_n(\mathbf{x}_n)) = \min_{D \in \hat{E}_n(\mathbf{x}_n)} \gamma(C, D).$$

How small should $\gamma^*$ be to say that breakdown of $C$ has occurred? The usual choice in robust statistics would be the worst possible value. In the present setup, this value depends on the dataset and on the clustering method. For example, if $n = 12$ and $|C| = 6$, the minimum for $\gamma^*(C, E_n^*(\mathbf{x}_{n+g}))$ is $1/4$. It is attained by building $s^* = 6$ clusters with two points each, one of which is in $C$. But $s < 6$ may be fixed, and this would result in a larger minimum. Even under estimated $s$ the minimum may be larger. For example, if the points lie on the real line and the clustering method produces only connected clusters, we get $\gamma^*(C, E_n^*(\mathbf{x}_{n+g})) \geq 2/7$. In general, the worst possible value may be difficult to compute and sometimes only attainable by tricky combinatorics, while one would judge a cluster as "broken down" already in much simpler constellations of $E_n^*(\mathbf{x}_{n+g})$. I propose

(3.2)   $\gamma^* \leq \frac{2}{3} = \gamma(\{x, y\}, \{x\}) = \gamma(C, C_1)$    if $C_1 \subset C$, $|C_1| = |C|/2$,

as the breakdown condition motivated by the following lemma, which means that under this condition every cluster can break down, at least in the absence of further subtle restrictions on the possible clusterings.

LEMMA 3.3. *Let $E_n(\mathbf{x}_n) \ni C$ be a partition with $|E_n(\mathbf{x}_n)| \geq 2$. Let $\mathcal{S} \subseteq \mathbb{N}$ be the set of possible cluster numbers containing at least one element $s \geq 2$. Let $\mathcal{F} = \{F \text{ partition on } \mathbf{x}_n : |F| \in \mathcal{S}\}$. Then $\exists \hat{F} \in \mathcal{F} : \gamma^*(C, \hat{F}) \leq 2/3$, where $2/3$ is the smallest value for this to hold.*

DEFINITION 3.4. Let $(E_n)_{n \in \mathbb{N}}$ be a sequence of GCM's. The *classification breakdown point* of a cluster $C \in E_n(\mathbf{x}_n)$ is defined as

$$B_n^c(E_n, \mathbf{x}_n, C) = \min_g \left\{ \frac{g}{n+g} : \exists \mathbf{x}_{n+g} = (x_1, \ldots, x_{n+g}) : \gamma^*(C, E_n^*(\mathbf{x}_{n+g})) \leq \frac{2}{3} \right\}.$$

The *$r$-clusters classification breakdown point* of $E_n$ at $\mathbf{x}_n$ is

$$B_n^c(E_n, \mathbf{x}_n) = \min_g \left\{ \frac{g}{n+g} : \exists \mathbf{x}_{n+g} = (x_1, \ldots, x_{n+g}), C_1, \ldots, C_r \in E_n(\mathbf{x}_n) \right.$$
$$\left. \text{pairwise distinct} : \gamma^*(C_i, E_n^*(\mathbf{x}_{n+g})) \leq \frac{2}{3}, i = 1, \ldots, r \right\}.$$

REMARK 3.5. At least $r \geq 1$ clusters of $E_n(\mathbf{x}_n)$ have to break down if $|E_n^*(\mathbf{x}_{n+g})| = s - r$. For $r = 1$, the reason is that there must be $D \in E_n^*(\mathbf{x}_{n+g})$ such that there are at least two members of $E_n(\mathbf{x}_n)$, $C_1$ and $C_2$, say, for which $D$ minimizes $\gamma(C_j, D)$ over $E_n^*(\mathbf{x}_{n+g})$. Without loss of generality, $|C_1 \cap D| \leq |C_2 \cap D|$. Since $C_1 \cap C_2 = \varnothing$, we get $\gamma(C_1, D) \leq |D|/(|D|/2 + |D|)$. The same argument yields for $r > 1$ that $\gamma(C_j, D) \leq 2/3$ for at least $q - 1$ clusters $C_j$ if $D \in E_n^*(\mathbf{x}_{n+g})$ is the most similar cluster for $q$ clusters $C_j \in E_n(\mathbf{x}_n)$.

Note that parameter breakdown does not imply classification breakdown and vice versa (cf. Remarks 4.10 and 4.18).

## 4. Breakdown results.

4.1. *Breakdown points for fixed $s$.* This section starts with three lemmas which characterize the behavior of the estimators for a sequence of datasets where there are $s \geq h \geq 2$ groups of points in every dataset, each group having a fixed range, but with the distances between the groups converging to $\infty$. In this case, eventually there exists a mixture component corresponding to each group, all mixture components correspond to one of the groups and the maximum of the log-likelihood can be obtained from the maxima considering the groups alone; that is, all groups are fitted separately.

LEMMA 4.1. *Let $\mathbf{x}_{nm} = (x_{1m}, \ldots, x_{nm}) \in \mathbb{R}^n$ be a sequence of datasets with $m \in \mathbb{N}$ and $0 = n_0 < n_1 < \cdots < n_h = n$, $h \geq 1$. Let $D_1 = \{1, \ldots, n_1\}$,*

$D_2 = \{n_1+1, \ldots, n_2\}, \ldots, D_h = \{n_{h-1}+1, \ldots, n_h\}$. *Assume further that*

$$\exists b < \infty : \max_k \max_{i,j \in D_k} |x_{im} - x_{jm}| \leq b \quad \forall m,$$

$$\lim_{m \to \infty} \min_{k \neq l, i \in D_k, j \in D_l} |x_{im} - x_{jm}| = \infty.$$

*Let* $s \geq h$ *be fixed,* $\eta_m = \arg\max_\eta L_{n,s}(\eta, \mathbf{x}_{nm})$. *The parameters of* $\eta_m$ *are called* $\pi_{1m}, \ldots, \pi_{sm}, a_{1m}$ *and so on; all results hold for* $\zeta_m$ *from maximizing* (2.10) *as well. Without loss of generality, assume* $x_{1m} \leq x_{2m} \leq \cdots \leq x_{nm}$. *Then, for* $m_0 \in \mathbb{N}$ *large enough,*

$$\exists 0 \leq d < \infty, \pi_{\min} > 0, \sigma_0 \leq \sigma_{\max} < \infty : \forall m > m_0,$$

(4.1) $\quad k = 1, \ldots, h \,\exists\, j_k \in \{1, \ldots, s\} : a_{j_k m} \in [x_{(n_{k-1}+1)m} - d, x_{n_k m} + d]$,

$$\pi_{j_k m} \geq \pi_{\min}, \sigma_{j_k m} \in [\sigma_0, \sigma_{\max}].$$

LEMMA 4.2. *In the situation of Lemma* 4.1, *assume further*

(4.2) $\quad \exists\, \pi_{\min} > 0 : \forall j = 1, \ldots, s, m \in \mathbb{N} : \pi_{jm} \geq \pi_{\min}.$

*Then*

(4.3) $\forall m > m_0, j = 1, \ldots, s \,\exists\, k \in \{1, \ldots, h\} : a_{jm} \in [x_{(n_{k-1}+1)m} - d, x_{n_k m} + d].$

(4.4) $\quad \exists 0 \leq \sigma_{\max} < \infty : \forall m > m_0, j = 1, \ldots, s : \sigma_{jm} \in [\sigma_0, \sigma_{\max}].$

LEMMA 4.3. *Under the assumptions of Lemma* 4.1,

(4.5) $\quad \forall k \in \{1, \ldots, h\} : \lim_{m \to \infty} \sum_{a_{jm} \in [x_{(n_{k-1}+1)m} - d, x_{n_k m} + d]} \pi_{jm} = \frac{|D_k|}{n},$

(4.6)
$$\lim_{m \to \infty} \left| L_{n,s}(\eta_m, \mathbf{x}_{nm}) - \max_{\sum_{k=1}^h q_k = s} \left( \sum_{k=1}^h \left[ \max_\eta L_{|D_k|, q_k}(\eta, \mathbf{y}_{km}) + |D_k| \log \frac{|D_k|}{n} \right] \right) \right| = 0,$$

*where* $\mathbf{y}_{km} = (x_{(n_{k-1}+1)m}, \ldots, x_{n_k m})$, $k = 1, \ldots, h$.

In particular, $r < s$ added outliers let $r$ mixture components break down if the differences between them tend to $\infty$.

THEOREM 4.4. *Let* $\mathbf{x}_n \in \mathbb{R}^n$, $s > 1$. *Let* $\eta_{n,s}$ *be a global maximizer of* (2.9). *Assume* (2.2)–(2.5). *For* $r = 1, \ldots, s-1$,

(4.7) $$B_{r,n}(\eta_{n,s}, \mathbf{x}_n) \leq \frac{r}{n+r}.$$

Equality in (4.7) could be proven for datasets where $\pi_j \to 0$ can be prevented for $j = 1, \ldots, s$ and any sequence of sets of $r$ added points, but conditions for this are hard to derive.

Under the restriction (2.12), convergence of $\sigma_j$-parameters to 0 implies breakdown according to Definition 3.1. Thus, to prevent breakdown, an effective lower bound for the $\sigma_j$ of the nonbreaking components has to exist. This means that *all* $\sigma_j$ have to be bounded from below, independently of $x_{n+1}, \ldots, x_{n+r}$, because (2.12) forces all $\sigma_j$ to 0 if only one implodes. Therefore, the result carries over.

COROLLARY 4.5. *Theorem 4.4 holds as well under the restriction* (2.12) *instead of* (2.11).

REMARK 4.6. The situation for $r = s$ is a bit more complicated, because here the choice of the basic distribution $f$ matters. Assume that the proportion of outliers in the distorted dataset is smaller than $1/2$. While $s - 1$ mixture components can be broken down by $s$ outliers, there remains at least one mixture component for which the original points own a majority of the weights used for the estimation of the parameters. If the parameters of such a component are estimated by a nonrobust ML-estimator such as the Normal one, the $s$th component will break down as well, that is, under $f = \varphi$,

$$B_{s,n}(\eta_{n,s}, \mathbf{x}_n) \leq \frac{s}{n+s}.$$

The breakdown point for the joint ML-estimator of location and scale for a single location–scale model based on the $t_\nu$-distribution was shown to be greater than or equal $1/(\nu + 1)$ by Tyler (1994), ignoring the possible breakdown of the scale to 0, which is prevented here because of (2.11). Suppose that points are added so that their proportion is smaller than $1/(\nu + 1)$. Mixture ML-estimation with $s$ components leads to the existence of at least one component such that the parameters are estimated by a weighted $t_\nu$-likelihood according to (2.17) with weight proportion smaller than $1/(\nu + 1)$ for the added points. Thus,

$$B_{s,n}(\eta_{n,s}, \mathbf{x}_n) \geq \frac{1}{\nu + 1},$$

where $f(x) = q(1 + x^2/\nu)^{(-\nu+1)/2}$, $\nu \geq 1$, $q > 0$ being the norming constant.

The approach via adding a noise component does not lead to a better breakdown behavior, because a single outlier can make the density value of the noise component arbitrarily small.

COROLLARY 4.7. *Theorem 4.4 and Remark 4.6 hold as well for global maximizers of* (2.10).

EXAMPLE 4.8. While the breakdown point for all considered approaches is the same for $r < s$, it may be of interest to determine how large an outlier has to be to cause breakdown of the methods. The following definition is used to generate reproducible datasets.

DEFINITION 4.9. $\Phi_{a,\sigma^2}^{-1}(1/(n+1)), \ldots, \Phi_{a,\sigma^2}^{-1}(n/(n+1))$ is called an $(a, \sigma^2)$-*Normal standard dataset* (*NSD*) with $n$ points, where $\Phi_{a,\sigma^2}$ denotes the c.d.f. of the Normal distribution with parameters $a$, $\sigma^2$.

Consider a dataset of 50 points, consisting of a $(0,1)$-NSD with 25 points combined with a $(5,1)$-NSD with 25 points (see Figure 1) and $s = 2$. For Normal mixtures, $t_\mu$-mixtures with $\mu \geq 1$ and Normal mixtures with noise component, the ML-estimators always result in components corresponding almost exactly to the two NSD's under $\sigma_0 = 0.025$ (see Section A.1). How large does an additional outlier have to be chosen so that the 50 original points fall into one single cluster and the second mixture component fits only the outlier? For Normal mixtures, breakdown begins with an additional point at about 15.2. For a mixture of $t_3$-distributions the outlier has to lie at about 800, $t_1$-mixtures need the outlier at about $3.8 \times 10^6$ and a Normal mixture with an additional noise component breaks down with an additional point at $3.5 \times 10^7$. These values depend on $\sigma_0$.

REMARK 4.10. Theorem 4.4 and Corollary 4.7 carry over to the classification breakdown point. This follows because if $r$ outliers are added, tending to $\infty$ and with the distance between them converging to $\infty$ as well, Lemma 4.2 yields that $p_{ij} \to 0$ for the original points $i = 1, \ldots, n$ and $j$ satisfying $a_{jm} \in [x_{n+g} - d, x_{n+g} + d]$ for some $g \in \{1, \ldots, r\}$. Thus, at most $s - r$ clusters remain for the classification of the original points, which yields breakdown of $r$ clusters; compare Remark 3.5. In contrast, the arguments leading to Remark 4.6 (Normal case) do not carry over because the addition of $r = s$ outliers as above certainly causes all mean parameters to explode, but one cluster usually remains containing all the original points. Therefore, an original cluster containing more than half of the points does not break down in the sense of classification.

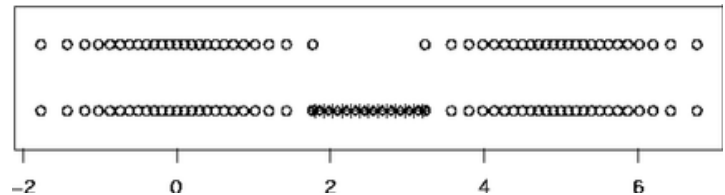

FIG. 1. Above: "*Standard*" *example dataset:* 25 *points* $(0,1)$-*NSD combined with* 25 *points* $(5,1)$-*NSD.* Below: *Stars denote* 13 *additional equidistant points between* 1.8 *and* 3.2.

4.2. *Alternatives for fixed s.* The results given above indicate that the considered mixture methods are generally not breakdown robust for fixed $s$. A first proposal for the construction of estimators with better breakdown behavior is based on the optimization of a target function for only a part of the data, say, optimally selected 50% or 80% of the points. The methods of trimmed $k$-means [Garcia-Escudero and Gordaliza (1999)] and clustering based on minimum covariance determinant estimators [Rocke and Woodruff (2000) and Gallegos (2003)] use this principle. Both methods, however, assume a partition model as opposed to the mixture model. Such an assumption may be useful for clustering, but yields biased parameter estimators [Bryant and Williamson (1986)]. Weighted likelihood as proposed by Markatou (2000) might be an alternative for the mixture model. One of the estimators treated in the previous section might be used after removing outliers by the nearest neighbor clutter removal procedure (NNC) of Byers and Raftery (1998). However, this procedure is based on mixture estimation as well (though not of location–scale type), and arguments analogous to those given above will lead to similar breakdown properties. As a simple example, consider a dataset consisting of a $(0, 1)$-NSD with 25 points, a $(5, 1)$-NSD with 25 points and an outlier at 50. The outlier at 50 is classified as "clutter" by NNC, but if another outlier as huge as $10^{100}$ is added, NNC classifies 50 as a nonoutlier.

Another alternative can be constructed by modifying the uniform noise approach. The problem of this approach is that the noise component could be affected by outliers as well, as was shown in the previous section. This can be prevented by choosing the density constant for the noise component as fixed in advance, leading to ML-estimation for a mixture where some improper distribution component is added to model the noise. That is, an estimator $\xi_{n,s}$ of the mean and variance parameters of the nonnoisy mixture components and of all the proportions is defined as the maximizer of

$$(4.8) \qquad L_{n,s}(\xi, \mathbf{x}_n) = \sum_{i=1}^{n} \log\left(\sum_{j=1}^{s} \pi_j f_{a_j, \sigma_j}(x_i) + \pi_0 b\right),$$

where $b > 0$. The choice of $b$ is discussed in Section A.1. For $\xi_{n,s}$, the breakdown point depends on the dataset $\mathbf{x}_n$. Breakdown can only occur if additional observations allow the nonoutliers to be fitted by fewer than $s$ components, and this means that a relatively good solution for $r < s$ components must exist even for $\mathbf{x}_n$. This is formalized in the following theorem, where only the breakdown of a single mixture component $B_{1,n}(\xi_{n,s}, \mathbf{x}_n)$ is considered.

THEOREM 4.11. *Let* $L_{n,s} = L_{n,s}(\xi_{n,s}, \mathbf{x}_n)$, $\mathbf{x}_n \in \mathbb{R}^n$. *Let* $\xi = \xi_{n,s}$ *and* $f_{\max} = f(0)/\sigma_0 > b$. *If*

$$
\max_{r<s} L_{n,r} < \sum_{i=1}^{n} \log\left(\sum_{j=1}^{s} \pi_j f_{\theta_j}(x_i) + \left(\pi_0 + \frac{g}{n}\right)b\right)
$$
(4.9)
$$
+ g \log\left(\pi_0 + \frac{g}{n}\right)b + (n+g)\log\frac{n}{n+g} - g\log f_{\max},
$$

*then*

(4.10) $$B_{1,n}(\xi_{n,s}, \mathbf{x}_n) > \frac{g}{n+g}.$$

EXAMPLE 4.12. Consider the dataset of 50 points shown in Figure 1, $f = \phi$, $b = 0.0117$ and $\sigma_0 = 0.025$ (cf. Section A.1). This results in $L_{n,1} = -119.7$. Neither the optimal solution for $s = 1$ nor the one for $s = 2$ classifies any point as noise. The right-hand side of (4.9) equals $-111.7$ for $g = 1$ and $-122.4$ for $g = 2$. Thus, the breakdown point is greater than $1/51$. Empirically, the addition of three extreme outliers at value 50, say, leads to a breakdown, namely to the classification of one of the two original components as noise and to the interpretation of the outliers as the second normal component. Two outliers do not suffice. Equation (4.10) is somewhat conservative. This stems from the exclusion of the breakdown of a proportion parameter to 0, which is irrelevant for this example.

A more stable data constellation with two clusters is obtained when a $(50,1)$-NSD of 25 points is added to the $(0,1)$-NSD of the same size. The optimal solution for one cluster classifies one of the two NSD's as noise and the other one as the only cluster, while the optimal solution for two clusters again does not classify any point as noise. Equation (4.9) leads to a minimal breakdown point of $8/58$ for the two-cluster solution. At least 11 outliers (at 500, say) are needed for empirical breakdown.

4.3. *Unknown s.* The treatment of the number of components $s$ as unknown is favorable for robustness against outliers, because outliers can be fitted by additional mixture components. Generally, for large enough outliers the addition of a new mixture component for each outlier yields a better log-likelihood than any essential change of the original mixture components. Thus, gross outliers are almost harmless, except that they let the estimated number of components grow.

Breakdown may occur, however, because additional points inside the range of the original data may lead to a a solution with $r < s$ clusters. Equation (4.11) of Theorem 4.13 is sufficient (but rather conservative) for preventing this. Breakdown can also occur due to gross outliers alone, simply because the number of outliers becomes so large that the BIC penalty, which depends on $n$, is increased by so much that the whole original dataset implodes into fewer than $s$ clusters. The conditions for this are given in (4.13) for BIC, while it cannot happen for AIC because its penalty does not depend on $n$.

THEOREM 4.13. *Let $\tau_n = (s, \eta_{n,s})$ be a maximizer of BIC. If*

$$(4.11) \quad \min_{r<s}[L_{n,s} - L_{n,r} - \tfrac{1}{2}(5g + 3s - 3r + 2n)\log(n+g) + n\log n] > 0,$$

*then*

$$(4.12) \qquad\qquad B_n(\tau_n, \mathbf{x}_n) > \frac{g}{n+g}.$$

*If*

$$(4.13) \qquad\qquad \min_{r<s}[L_{n,s} - L_{n,r} - \tfrac{3}{2}(s-r)\log(n+g)] < 0,$$

*then*

$$(4.14) \qquad\qquad B_n(\tau_n, \mathbf{x}_n) \leq \frac{g}{n+g}.$$

Note that $L_{n,s} - L_{n,r} > 3/2(s-r)\log n$ always holds by definition of BIC. Sufficient conditions for breakdown because of "inliers" depend on the parameters of certain suboptimal solutions for $r \leq s$ mixture components for $\mathbf{x}_n$. They may be hard to derive and are presumably too complicated to be of practical use.

EXAMPLE 4.14. Consider again the combination of a $(0,1)$-NSD with 25 points and a $(5,1)$-NSD with 25 points, $f = \varphi$ and $\sigma_0$ chosen as in Example 4.12. The difference in (4.11) is 3.37 for $g = 1$ and $-7.56$ for $g = 2$; that is, the breakdown point is larger than $1/51$. Many more points are empirically needed. Thirteen additional points, equally spaced between 1.8 and 3.2, lead to a final estimation of only one mixture component (compare Figure 1). It may be possible to find a constellation with fewer points where one component fits better than two or more components, but I did not find any. Breakdown because of gross outliers according to (4.13) needs more than 650,000 additional points!

A mixture of the $(0,1)$-NSD with 25 points with a $(50,1)$-NSD of size 25 leads to a lower breakdown bound of $12/62$. For estimated $s$, even a breakdown point larger than $1/2$ is possible, because new mixture components can be opened for additional points. This may even happen empirically for a mixture of $(0,1)$-NSD and $(50,1)$-NSD, because breakdown by addition of gross outliers is impossible unless their number is huge, and breakdown by addition of "inliers" is difficult. For a $(0, 0.001)$-NSD of 25 points and a $(100{,}000, 0.001)$-NSD of 25 points, even the conservative lower breakdown bound is $58/108 > 1/2$.

The choice of the $t_1$-distribution instead of the Normal leads to slightly better breakdown behavior. The mixture of a 25 point-$(0,1)$-NSD and a 25 point-$(5,1)$-NSD yields a lower breakdown bound of $3/53$, and empirically the addition of the 13 inliers mentioned above does not lead to breakdown of one of the two components, but to the choice of three mixture components by the BIC. Replacement of the $(5,1)$-NSD by a $(50,1)$-NSD again gives a small improvement of the lower bound to $13/63$.

REMARK 4.15. The possible breakdown point larger than $1/2$ is a consequence of using the addition breakdown definition. A properly defined replacement breakdown point can never be larger than the portion of points in the smallest cluster, because this cluster must be driven to break down if all of its points are suitably replaced. This illustrates that the correspondence between addition and replacement breakdown as established by Zuo (2001) may fail in more complicated setups.

The addition of a noise component again does not change the breakdown behavior.

THEOREM 4.16. *Under $f_{\max} \geq 1/(x_{\max,n} - x_{\min,n})$, Theorem 4.13 also holds for global maximizers of BIC, defined so that (2.10) is maximized for every fixed $s$.*

EXAMPLE 4.17. The discussed data examples of two components with 25 points each do not lead to different empirical breakdown behavior with and without an estimated noise component according to (2.10), because no point of the original mixture components is classified as noise by the solutions for two Normal components. In the case of a $(0,1)$-NSD of 45 points and a $(5,1)$-NSD of 5 points, the solution with one Normal component, classifying the points from the smaller NSD as noise, is better than any solution with two components. That is, no second mixture component exists which could break down. The same holds for $t_1$-mixtures (all points form the only component), while NMML shows almost the same behavior in Example 4.14: there are two mixture components corresponding to the two NSD's which can be joined by 12 equidistant points between 1.55 and 3.55. Equation (4.12) evaluates again to $1/51$. More examples are given in Hennig (2003).

REMARK 4.18. While parameter breakdown due to the loss of a mixture component implies classification breakdown of at least one cluster, classification breakdown may occur with fewer additional points than parameter breakdown. Consider again the $(0,1)$-NSD of 45 points plus the $(5,1)$-NSD of 5 points and NMML. The smaller cluster breaks down by the addition of six points, namely two points each exactly at the smallest and the two largest points of the $(5,1)$-NSD. This leads to the estimation of five clusters, namely the original $(0,1)$-NSD, three clusters of three identical points each, and the remaining two points of the $(5,1)$-NSD. The fifth cluster is most similar to the original one with $\gamma = \frac{2*2}{2+5} < \frac{2}{3}$, while no parameter breakdown occurs. Thus, an arbitrarily large classification breakdown point is not possible even for very well separated clusters, because not only their separation, but also their size matters. As in Section 4.2, the number of additional points required depends on $\sigma_0$.

**5. Discussion.** It has been shown that none of the discussed mixture model estimators is breakdown robust when the number of components $s$ is assumed as known and fixed. An improvement can be achieved by adding an improper uniform distribution as an additional mixture component.

The more robust way of estimating mixture parameters is the simultaneous estimation of the number of mixture components $s$. Breakdown of mixture components may rather arise from the addition of points between the estimated mixture components of the original dataset than from gross outliers. It may be controversial if this is really a robustness problem. A sensible clustering method should be expected to reduce the estimated number of clusters if the gap between the clusters is filled with points, as long as their number is not too small. Compare Figure 1, where the NMML estimate of $s=1$ and the $t_1$-mixture estimate of $s=3$ may both seem to be acceptable. In such cases, the empirical breakdown point, or the more easily computable but conservative breakdown bound (4.12), may not be used to rule out one of the methods, but can rather be interpreted as a measure of the stability of the dataset with respect to clustering.

While including the estimation of $s$ leads to theoretically satisfying breakdown behavior, robustness problems remain, in practice, because the global optimum of the log-likelihood has to be found. Consider, for example, a dataset of 1000 points, consisting of three well-separated clusters of 300 points each and 100 extremely scattered outliers. The best solution requires 103 clusters. Even for one-dimensional data, however, the EM-algorithm will be very slow for a large number of clusters, and there will be typically lots of local optima. Therefore, the maximum number of fitted components will often be much smaller than the maximum possible number of outliers and the results for fixed $s$ remain relevant. The use of an improper noise component or, if extremely huge outliers are ruled out, the proper noise component or $t_1$-mixtures will clearly be superior to Normal mixtures with $s$ estimated but restricted to be small.

The comparison of the robustness characteristics of various cluster analysis methods is an important topic and a first attempt is made this paper to define a classification breakdown point. It should not be the last word on the subject. Davies and Gather (2002) argue that a reasonable concept of a breakdown point should be linked to a sufficiently rich equivariance structure to enable nontrivial upper bounds for the breakdown point. This does not hold for the concepts presented here, and it should be kept in mind that breakdown point definitions as those given here do not rule out meaningless estimators such as constants. The breakdown point should not be seen as the only important measure to judge the methods, but must be complemented by the consideration of their further properties.

In some situations with low breakdown point in mixture modeling, additional outliers do not cause any substantial change unless they are huge (cf. Example 4.8, NNC in Section 4.2). More sensible measures than the breakdown point may be needed here.

Neither MCLUST nor EMMIX is able to exactly reproduce the results given here. Both do not allow the specification of a lower scale bound. MCLUST produces an error if the EM-iteration leads to a sequence of variance parameters converging to 0. This implies, in particular, that no single point can be isolated as its own mixture component. But such an isolation is crucial for the desirable breakdown behavior of the methods with estimated $s$. EMMIX terminates the iteration when the log-likelihood does not seem to converge. The preliminary iteration results, including one-point-components, are reported, but solutions with clear positive variances are favored. Thus, the current implementations of the Normal mixture estimation with estimated $s$ are essentially nonrobust. Addition of a noise component and $t$-mixtures perform better under outliers of moderate size, but they, too, are not robust against very extreme outliers. The results given here do not favor one of these two approaches over the other, and I think that the implementation of a lower bound for the smallest covariance eigenvalue is more important an issue than the decision between the current implementations.

Note that both packages enable the use of stronger scale restrictions (equivalent to equal variances for all mixture components in the one-dimensional case), which should have roughly the same robustness characteristics for estimated $s$ as the methods considered here. However, in practice such restrictions are often not justified.

## APPENDIX

**A.1. Choice of the tuning parameters $\sigma_0$ and $b$.** For the choice of $\sigma_0$, the following strategy is proposed. As a "calibration benchmark," form a dataset with $n$ points by adding an $\alpha_n$-outlier to a $(0,1)$-NSD (recall Definition 4.9) with $n-1$ points. Davies and Gather (1993) define "$\alpha$-outliers" (with $\alpha > 0$ but very small) with respect to an underlying model as points from a region of low density, chosen so that the probability of the occurrence of an outlier is equal to $\alpha$ under that model. For a standard Normal distribution, for example, the points outside $[\Phi^{-1}(\frac{\alpha}{2}), \Phi^{-1}(1 - \frac{\alpha}{2})]$ are the $\alpha$-outliers. For $\alpha_n = 1 - (1-p)^{1/n}$, the probability of the occurrence of at least one $\alpha_n$-outlier among $n$ i.i.d. points is equal to $p$. Take $p = 0.95$, say.

Consider NMML with estimated $s$ under (2.11) (this seems to lead to reasonable values for all methods discussed in the present paper). Let $c_0 = \sigma_0$ for this particular setup. Choose $c_0$ so that $C(1) = C(2)$ according to (2.20). This can be carried out in a unique way because $L_{n,1}(\eta_{n,1})$ does not depend on $c_0$ (as long as $c_0$ is smaller than the sample variance) and $L_{n,2}(\eta_{n,2})$ increases with decreasing $c_0$, because this enlarges the parameter space. For $c_0$ small enough, the two-component solution will consist of one component matching approximately the ML-estimator for the NSD, $a_2$ will approximately equal the outlier and $\sigma_2 = c_0$, so that the increase in $L_{n,2}(\eta_{n,2})$ becomes strict.

Now use $\sigma_0 = c_0 \sigma_{\max}$, where $\sigma_{\max}^2$ is the largest variance such that a data subset with this variance can be considered as a "cluster" with respect to the given application. At least, if the mixture model is used as a tool for cluster analysis, points of a cluster should belong together in some sense, and, with regard to a particular application, it can usually be said that points above a certain variation can no longer be considered as "belonging together." Therefore, in most applications it is possible to choose $\sigma_{\max}$ in an interpretable manner, while this does not work for $\sigma_0$ directly.

The rationale is that a sensible choice of $\sigma_0$ should lead to the estimation of the dataset as one component, if it does not contain any outlier in the sense of Davies and Gather (1993). If the $n$th point is an outlier, it should be fitted by a new mixture component. The reader is referred to Hennig (2003) for a more detailed discussion.

Given $\sigma_{\max}$, the improper density value $b$ for maximization of (4.8) can be chosen as the density value at the 0.025-quantile of $f_{0,\sigma_{\max}}$, so that at least 95% of the points generated from a "cluster-generating" mixture component have a larger density value for their own parent distribution than for the noise component. In all examples $\sigma_{\max} = 5$ has been used, which leads to $\sigma_0 = 0.025$, $b = 0.0117$.

Note that the theory in Section 4 assumes $\sigma_0$ as constant over $n$, so that it does not directly apply to the suggestion given here.

Under (2.12), $c = c_0$ can be used because of scale equivariance, avoiding the specification of $\sigma_{\max}$. However, (2.12) does not properly generalize to fitting of a noise component and estimation of the number of components (the latter can be done by the choice of a suitable upper bound on $s$).

LEMMA A.1. *The following objective functions are unbounded from above under the restriction* (2.12):

1. *the log-likelihood function* (2.10) *with fixed s;*
2. *the AIC and BIC of model* (2.1) *with unknown* $s \in \mathbb{N}$.

PROOF. Consider an arbitrary dataset $x_1, \ldots, x_n$. For (2.10) choose $a_1 = x_1$, $\pi_1 > 0$, $\sigma_1 \to 0$, $\pi_0 > 0$. This implies that the summand for $x_1$ converges to $\infty$ while all others are bounded from below by $\log(\pi_0/(x_{\max,n} - x_{\min,n}))$. This proves part (1). For part (2) choose $s = n$, $a_1 = x_1, \ldots, a_s = x_n$, $\sigma_1 = \cdots = \sigma_s \to 0$. Thus, $L_{n,s} \to \infty$, and the same holds for AIC and BIC. □

### A.2. Proofs.

PROOF OF LEMMA 2.2. For any fixed $\sigma_j^*$, the maximizer $a_j$ of (2.17) lies between $x_{\max,n}$ and $x_{\min,n}$ because of (2.2) and (2.3). Now show that

$$\sigma_j \leq \frac{\sigma_0 f(0)}{f((x_{\max,n} - x_{\min,n})/\sigma_0)}.$$

By $\sigma_j^* = \sigma_0$,

$$S_j(a_j, \sigma_j) \geq \sum_{i=1}^n p_{ij} \log \frac{1}{\sigma_0} f\left(\frac{x_i - a_j}{\sigma_0}\right) \geq n\pi_j \log \frac{1}{\sigma_0} f\left(\frac{x_{\max,n} - x_{\min,n}}{\sigma_0}\right).$$

For arbitrary $\sigma_j^*$,

$$S_j(a_j, \sigma_j^*) \leq n\pi_j(\log f(0) - \log \sigma_j^*).$$

Therefore,

$$\log f(0) - \log \sigma_j \geq \log \frac{1}{\sigma_0} f\left(\frac{x_{\max,n} - x_{\min,n}}{\sigma_0}\right) \Rightarrow \sigma_j \leq \frac{\sigma_0 f(0)}{f((x_{\max,n} - x_{\min,n})/\sigma_0)}$$

as long as $n\pi_j > 0$, proving (2.18). □

PROOF OF LEMMA 3.3. Recall (3.2). For given $C$, $\hat{F}$ can always be chosen to contain $C_1, \ldots, C_r$, $r \geq 2$, with $C \subseteq \bigcup_{i=1}^r C_i$ such that $|C_i| \leq |C|/2$ $\forall i$ for even $|C|$. For odd $|C|$, $C_1 \in \hat{F}$ with $|C_1 \cap C| = (|C|+1)/2$ and $|C_1 \setminus C| \geq 1$ can be constructed such that $\gamma^*(C, \hat{F}) = \gamma(C, C_1) \leq 2/3$. On the other hand, $\forall \hat{F} \in \mathcal{F}: \gamma^*(C, \hat{F}) \geq 2/3$ if $\mathbf{x}_n = C \cup \{x\}$ and $\mathcal{S} = \{2\}$. □

PROOF OF LEMMA 4.1. Note first that in case of maximizing (2.10) the density of the noise component $1/(x_{\max,n+g} - x_{\min,n+g})$ converges to 0, so that all arguments, including those used in the proofs of Lemmas 4.2 and 4.3, hold for this case too.

Assume w.l.o.g. that all $a_{jm}$, $j = 1, \ldots, s$, are outside $[x_1 - d, x_{n_1} + d]$ for arbitrary $d < \infty$ and $m$ large enough unless $\pi_{jm} \searrow 0$ or $\sigma_{jm} \nearrow \infty$ at least for a subsequence of $m \in \mathbb{N}$. Consider

$$L_{n,s}(\eta_m, \mathbf{x}_{nm}) = \sum_{i=1}^{n_1} \log\left(\sum_{j=1}^{s} \pi_{jm} f_{a_{jm},\sigma_{jm}}(x_i)\right)$$
$$+ \sum_{i=n_1+1}^{n} \log\left(\sum_{j=1}^{s} \pi_{jm} f_{a_{jm},\sigma_{jm}}(x_i)\right).$$

The first sum converges to $-\infty$ for $m \to \infty$ because of (2.7), and the second sum is bounded from above by $(n-n_1)\log(f(0)/\sigma_0)$, that is, $L_{n,s}(\eta_\mu, \mathbf{x}_{nm}) \to -\infty$. In contrast, for $\hat\eta_m$ with $\hat{a}_{km} = x_{n_k}$, $\hat\sigma_{km} = \sigma_0$, $\hat\pi_{km} = \frac{1}{h}$, $k = 1, \ldots, h$,

$$L_{n,s}(\hat\eta_m, \mathbf{x}_{nm}) \geq \sum_{k=1}^{h} n_k \log \frac{f((x_{n_k m} - x_{(n_{k-1}+1)m})/\sigma_0)}{h\sigma_0} \geq n \log \frac{f(b/\sigma_0)}{h\sigma_0} > -\infty.$$

Hence, for $m$ large enough, $\eta_m$ cannot be ML. Since it should be ML, $d$ has to exist so that (4.1) holds for $m$ larger than some $m_0$. $\square$

PROOF OF LEMMA 4.2.

*Proof of* (4.3). Suppose that (4.3) does not hold. Without loss of generality [the order of the $a_j$ does not matter and a suitable subsequence of $(\eta_m)_{m \in \mathbb{N}}$ can always be found] assume

$$\lim_{m \to \infty} \min\{|x - a_{1m}| : x \in \{x_{1m}, \ldots, x_{nm}\}\} = \infty.$$

Due to (2.7),

$$\frac{1}{\sigma_{1m}} f\left(\frac{x_{im} - a_{1m}}{\sigma_{1m}}\right) \to 0 \quad \forall i.$$

With (2.6) and (4.1),

$$\sum_{j=2}^{s} \pi_{jm} f_{a_{jm},\sigma_{jm}}(x_i) \geq d_{\min} = \pi_{\min} \frac{1}{\sigma_{\max}} f\left(\frac{b+2d}{\sigma_0}\right) > 0, \quad i = 1, \ldots, n.$$

Thus, for arbitrarily small $\varepsilon > 0$ and $m$ large enough,

$$L_{n,s}(\eta_m, \mathbf{x}_{nm}) \leq \sum_{i=1}^{n} \log\left(\sum_{j=2}^{s} \pi_{jm} f_{a_{jm},\sigma_{jm}}(x_i)\right) + n(\log(d_{\min} + \varepsilon) - \log d_{\min}),$$

and $\log(d_{\min} + \varepsilon) - \log d_{\min} \searrow 0$ for $\varepsilon \searrow 0$. Thus, $L_{n,s}$ can be increased for $\varepsilon$ small enough by replacement of $(\pi_{1m}, a_{1m}, \sigma_{1m})$ by $(\pi_{1m}, x_1, \sigma_0)$ in contradiction to $\eta_m$ being ML.

*Proof of* (4.4) *by analogy to* (4.3). Suppose that w.l.o.g. $\sigma_{1m} \to \infty$. Then
$$\frac{1}{\sigma_{1m}} f\left(\frac{x_{im} - a_{1m}}{\sigma_{1m}}\right) \to 0 \qquad \forall i,$$
and replacement of $(\pi_{1m}, a_{1m}, \sigma_{1m})$ by $(\pi_{1m}, x_1, \sigma_0)$ increases the log-likelihood. □

PROOF OF LEMMA 4.3.

*Proof of* (4.5). Consider $k \in \{1, \ldots, h\}$. Let $S_k = [x_{(n_{k-1}+1)m} - d, x_{n_k m} + d]$. With Lemma 2.2,
$$\sum_{a_{jm} \in S_k} \pi_{jm} = \sum_{a_{jm} \in S_k} \frac{1}{n} \sum_{i=1}^{n} \frac{\pi_{jm} f_{a_{jm}, \sigma_{jm}}(x_i)}{\sum_{l=1}^{s} \pi_{lm} f_{a_{lm}, \sigma_{lm}}(x_i)}.$$

For $a_{jm} \in S_k$ and $m \to \infty$,
$$\frac{\pi_{jm} f_{a_{jm}, \sigma_{jm}}(x_i)}{\sum_{l=1}^{s} \pi_{lm} f_{a_{lm}, \sigma_{lm}}(x_i)} \to 0 \qquad \text{for } i \notin D_k,$$
while, for $i \in D_k$,
$$\left| \frac{\pi_{jm} f_{a_{jm}, \sigma_{jm}}(x_i)}{\sum_{l=1}^{s} \pi_{lm} f_{a_{lm}, \sigma_{lm}}(x_i)} - \frac{\pi_{jm} f_{a_{jm}, \sigma_{jm}}(x_i)}{\sum_{a_{jm} \in S_k} \pi_{lm} f_{a_{lm}, \sigma_{lm}}(x_i)} \right| \to 0.$$

This yields $\sum_{a_{jm} \in S_k} \pi_{jm} \to |D_k|/n$ [at least one of the $\pi_{jm}$ in this sum is bounded away from 0 by (4.1)].

*Proof of* (4.6). Let $\eta_{kmq} = \arg\max_{\eta} L_{|D_k|, q}(\eta, \mathbf{y}_{km}), q \in \mathbb{N}$,
$$L_{q_1 \cdots q_h m} = \sum_{k=1}^{h} \left( L_{|D_k|, q_k}(\eta_{kmq_k}) + |D_k| \log \frac{|D_k|}{n} \right).$$

Note that $L_{n,s}(\eta_m) \geq \max_{\sum_{k=1}^{h} q_k = s} L_{q_1 \cdots q_h m}$ can be proved by choice of $\eta$ according to $\pi_j = (|D_j|/n) \pi_{jmq_j}$, $a_j = a_{jmq_j}$, $\sigma_j = \sigma_{jmq_j}$, $j = 1, \ldots, h$. Further, for $m$ large enough and arbitrarily small $\varepsilon > 0$,

$$(\text{A.1}) \qquad L_{n,s}(\eta_m) \leq \sum_{k=1}^{h} \sum_{i \in D_k} \log\left(\sum_{a_{jm} \in S_k} \pi_{jm} f_{a_{jm}, \sigma_{jm}}(x_i)\right) + \varepsilon,$$

because, for $x_i, i \in D_k$, the sum over $a_{jm} \in S_k$ is bounded away from 0 as shown in the proof of Lemma 4.1, while the sum over $a_{jm} \in [x_{(n_{l-1}+1)m} - d, x_{n_l m} + d]$, $l \neq k$, vanishes for $m \to \infty$. Further, find

$$\sum_{i \in D_k} \log\left(\sum_{a_{jm} \in S_k} \pi_{jm} f_{a_{jm}, \sigma_{jm}}(x_i)\right) - |D_k| \log\left(\sum_{a_{jm} \in S_k} \pi_{jm}\right) \leq L_{|D_k|, q}(\eta_{kmq}),$$

where $q = |\{a_{jm} \in S_k\}|$. Now (4.6) follows from (4.5). □

PROOF OF THEOREM 4.4. Let $\mathbf{x}_{(n+r)m} = (x_1, \ldots, x_n, x_{(n+1)m}, \ldots, x_{(n+r)m})$, $m \in \mathbb{N}$ w.l.o.g. Let $x_1 \leq \cdots \leq x_n$, $x_{(n+k)m} = x_n + km$, $k = 1, \ldots, r$. This satisfies the assumptions of Lemma 4.1 for $h = r+1$, so that the location parameters for $r$ components have to converge to $\infty$ with $x_{(n+1)m}, \ldots, x_{(n+r)m}$. □

PROOF OF THEOREM 4.11. Let $\mathbf{x}_{n+g} = (x_1, \ldots, x_{n+g})$. Let $\xi^* = \xi_{n+g,s} = \arg\max_{\hat{\xi}} L_{n+g,s}(\hat{\xi}, \mathbf{x}_{n+g})$. For $r < s$,

$$L_{n+g,s} \leq \sum_{i=1}^n \log\left(\sum_{j=1}^r \pi_j^* f_{\theta_j^*}(x_i) + \sum_{j=r+1}^s \pi_j^* f_{\theta_j^*}(x_i) + \pi_0^* b\right) + g \log f_{\max}.$$

Assume that the parameter estimators of $s - r$ (i.e., at least one) mixture components leave a compact set $D$ of the form $D = [\pi_{\min}, 1] \times C$, $C \subset \mathbb{R} \times \mathbb{R}^+$ compact, $\pi_{\min} > 0$. Let the mixture components be ordered such that $(\pi_j^*, a_j^*, \sigma_j^*) \in D$ only for $j = 1, \ldots, r < s$. From (2.6), $\sum_{j=1}^r \pi_j^* f_{\theta_j^*}(x_i) \geq r \pi_{\min} f_{\min}$, while $\sum_{j=r+1}^s \pi_j^* f_{\theta_j^*}(x_i)$ becomes arbitrarily small for $D$ large enough by (2.7). Thus, for arbitrary $\varepsilon > 0$ and $D$ large enough,

$$L_{n+g,s} \leq \sum_{i=1}^n \log\left(\sum_{j=1}^r \pi_j^* f_{\theta_j^*}(x_i) + \pi_0^* b\right) + g \log f_{\max} + \varepsilon$$

(A.2)
$$\leq \max_{r<s} L_{n,r} + g \log f_{\max} + \varepsilon.$$

However, $\hat{\xi}$ could be defined by $\hat{\pi}_0 = (n\pi_0 + g)/(n+g)$, $\hat{\pi}_j = (n/(n+g))\pi_j$, $\hat{a}_j = a_j$, $\hat{\sigma}_j = \sigma_j$, $j = 1, \ldots, s$. Therefore,

$$L_{n+g,s} \geq \sum_{i=1}^n \log\left(\sum_{j=1}^s \pi_j f_{\theta_j}(x_i) + \left(\pi_0 + \frac{g}{n}\right)b\right)$$

$$+ g \log\left[\left(\pi_0 + \frac{g}{n}\right)b\right] + (n+g) \log \frac{n}{n+g}$$

$$\Rightarrow \max_{r<s} L_{n,r} \geq \sum_{i=1}^n \log\left(\sum_{j=1}^s \pi_j f_{\theta_j}(x_i) + \left(\pi_0 + \frac{g}{n}\right)b\right)$$

$$+ g \log\left[\left(\pi_0 + \frac{g}{n}\right)b\right] + (n+g) \log \frac{n}{n+g} - g \log f_{\max} - \varepsilon.$$

This contradicts (4.9) by $\varepsilon \to 0$. □

PROOF OF THEOREM 4.13. Add points $x_{n+1}, \ldots, x_{n+g}$ to $\mathbf{x}_n$. Let $C_m(s, \hat{\eta})$ be the value of BIC for $s$ mixture components and parameter $\hat{\eta}$, applied to the dataset $\mathbf{x}_m$, $m \geq n$. Let $C_m(s)$ be its maximum. With the same arguments as those leading to (A.2), construct for arbitrary $\varepsilon > 0$ a suitably large compact $C \subset \mathbb{R} \times \mathbb{R}^+$, containing the location and scale parameters of all mixture components of $\tau = (s, \eta) = (s, \eta_{n,s})$, and assume that $(a_j^*, \sigma_j^*) \in C$ for only $r < s$ components of $\tau^* = \arg\max_{\hat{s}, \hat{\eta}} C_{n+g}(\hat{s}, \hat{\eta})$. We get

$$
\text{(A.3)} \quad
\begin{aligned}
C_{n+g}(s^*) \leq {}& 2 \sum_{i=1}^{n} \log\left( \sum_{j=1}^{r} \pi_j^* f_{\theta_j^*}(x_i) \right) \\
& + 2g \log f_{\max} + \varepsilon - (3s^* - 1) \log(n+g),
\end{aligned}
$$

and, by taking $\hat{s} = s + g$, $\hat{\pi}_j = n/(n+g)\pi_j$, $j = 1, \ldots, s$, $\hat{\pi}_{s+1} = \cdots = \hat{\pi}_{s+g} = 1/(n+g)$, $\hat{\theta}_j = \theta_j$, $j = 1, \ldots, s$, $\hat{a}_{s+k} = x_{n+k}$, $\hat{\sigma}_{s+k} = \sigma_0$, $k = 1, \ldots, g$,

$$
\text{(A.4)} \quad
\begin{aligned}
C_{n+g}(s^*) \geq {}& 2 \sum_{i=1}^{n} \log\left( \sum_{j=1}^{s} \frac{n}{n+g} \pi_j f_{\theta_j}(x_i) \right) \\
& + 2g \log \frac{f_{\max}}{n+g} - (3(s+g) - 1) \log(n+g).
\end{aligned}
$$

By combination,

$$
\sum_{i=1}^{n} \log\left( \sum_{j=1}^{s} \pi_j f_{\theta_j}(x_i) \right) - \sum_{i=1}^{n} \log\left( \sum_{j=1}^{r} \frac{\pi_j^*}{\sum_{k=1}^{r} \pi_k^*} f_{\theta_j^*}(x_i) \right) - \varepsilon
$$

$$
\leq g \log(n+g) - \frac{3}{2}(s^* - (s+g)) \log(n+g)
$$

$$
- n \log \frac{n}{n+g} + n \log\left( \sum_{k=1}^{r} \pi_k^* \right)
$$

$$
\leq \frac{1}{2}(5g + 3s - 3r + 2n) \log(n+g) - n \log n.
$$

Under (4.11) this cannot happen for arbitrarily small $\varepsilon$.

A sufficient condition for breakdown can be derived by explicit contamination. Let $y = x_{n+1} = \cdots = x_{n+g}$. For fixed $\hat{s}$, it follows from Lemma 4.3 that

$$
\lim_{y \to \infty} C_{n+g}(\hat{s}) = 2\left( L_{n, \hat{s}-1} + g \log(f_{\max}) + n \log \frac{n}{n+g} + g \log \frac{g}{n+g} \right)
$$
$$
- (3\hat{s} - 1) \log(n+g).
$$

This cannot be maximized by $s^* = \hat{s} > s + 1$ because the penalty on $s$ is larger for $n + g$ points than for $n$ points and $s^* - 1$ with parameters maximizing $L_{n, s^*-1}(\hat{\eta}, \mathbf{x}_n)$ must already be a better choice than $s$ for $n$ points unless $s^* \leq s + 1$. It follows that the existence of $r < s$ with

$$
2L_{n,s} - (3(s+1) - 1) \log(n+g) < 2L_{n,r} - (3(r+1) - 1) \log(n+g)
$$

suffices for breakdown of at least one component, which is equivalent to (4.13). □

PROOF OF THEOREM 4.16. Let
$$d = \frac{1}{x_{\max,n} - x_{\min,n}}, \qquad d^* = \frac{1}{x_{\max,n+g} - x_{\min,n+g}}.$$

Replace (A.3) by

$$(A.5) \quad C_{n+g}(s^*) \leq 2 \sum_{i=1}^{n} \log \left( \sum_{j=1}^{r} \pi_j^* f_{\theta_j^*}(x_i) + \pi_0^* d^* \right) + 2g \log f_{\max} + \varepsilon - (3s^* - 1)\log(n+g),$$

and (A.4) by

$$C_{n+g}(s^*) \geq 2 \sum_{i=1}^{n} \log \left( \sum_{j=1}^{s} \frac{n}{n+g} \pi_j f_{\theta_j}(x_i) + \frac{n}{n+g} \pi_0 d \right) + 2g \log \frac{f_{\max}}{n+g} - (3(s+g) - 1)\log(n+g).$$

Equation (4.12) follows from $d \geq d^*$ in (A.5). □

Lemma 4.3 holds as well for maximizers of (2.10), and therefore (4.14) carries over as well.

**Acknowledgments.** I thank a referee and Vanessa Didelez for helpful comments.

## REFERENCES

AKAIKE, H. (1974). A new look at the statistical model identification. *IEEE Trans. Automatic Control* **19** 716–723. MR423716

BANFIELD, J. D. and RAFTERY, A. E. (1993). Model-based Gaussian and non-Gaussian clustering. *Biometrics* **49** 803–821. MR1243494

BOZDOGAN, H. (1994). Mixture model cluster analysis using model selection criteria and a new informational measure of complexity. In *Multivariate Statistical Modeling. Proc. First US/Japan Conference on the Frontiers of Statistical Modeling. An Informational Approach* (H. Bozdogan, ed.) **2** 69–113. Kluwer, Dordrecht. MR1328880

BRYANT, P. and WILLIAMSON, J. A. (1986). Maximum likelihood and classification: A comparison of three approaches. In *Classification as a Tool of Research* (W. Gaul and M. Schader, eds.) 35–45. North-Holland, Amsterdam. MR913117

BYERS, S. and RAFTERY, A. E. (1998). Nearest neighbor clutter removal for estimating features in spatial point processes. *J. Amer. Statist. Assoc.* **93** 577–584.

CAMPBELL, N. A. (1984). Mixture models and atypical values. *Math. Geol.* **16** 465–477.

Celeux, G. and Soromenho, G. (1996). An entropy criterion for assessing the number of clusters in a mixture model. *J. Classification* **13** 195–212. MR1421665

Davies, P. L. and Gather, U. (1993). The identification of multiple outliers (with discussion). *J. Amer. Statist. Assoc.* **88** 782–801. MR1242933

Davies, P. L. and Gather, U. (2002). Breakdown and groups. Technical Report 57-2002, SFB 475, Univ. Dortmund. Available at wwwstat.mathematik.uni-essen.de/~davies/brkdown220902.ps.gz.

Day, N. E. (1969). Estimating the components of a mixture of normal distributions. *Biometrika* **56** 463–474. MR254956

Dempster, A. P., Laird, N. M. and Rubin, D. B. (1977). Maximum likelihood from incomplete data via the EM algorithm (with discussion). *J. Roy. Statist. Soc. Ser. B* **39** 1–38. MR501537

DeSarbo, W. S. and Cron, W. L. (1988). A maximum likelihood methodology for clusterwise linear regression. *J. Classification* **5** 249–282. MR971156

Donoho, D. L. and Huber, P. J. (1983). The notion of breakdown point. In *A Festschrift for Erich L. Lehmann* (P. J. Bickel, K. Doksum and J. L. Hodges, Jr., eds.) 157–184. Wadsworth, Belmont, CA. MR689745

Fraley, C. and Raftery, A. E. (1998). How many clusters? Which clustering method? Answers via model-based cluster analysis. *The Computer J.* **41** 578–588.

Gallegos, M. T. (2003). Clustering in the presence of outliers. In *Exploratory Data Analysis in Empirical Research* (M. Schwaiger and O. Opitz, eds.) 58–66. Springer, Berlin. MR2074222

Garcia-Escudero, L. A. and Gordaliza, A. (1999). Robustness properties of $k$ means and trimmed $k$ means. *J. Amer. Statist. Assoc.* **94** 956–969. MR1723291

Hampel, F. R. (1971). A general qualitative definition of robustness. *Ann. Math. Statist.* **42** 1887–1896. MR301858

Hampel, F. R. (1974). The influence curve and its role in robust estimation. *J. Amer. Statist. Assoc.* **69** 383–393. MR362657

Hastie, T. and Tibshirani, R. (1996). Discriminant analysis by Gaussian mixtures. *J. Roy. Statist. Soc. Ser. B* **58** 155–176. MR1379236

Hathaway, R. J. (1985). A constrained formulation of maximum-likelihood estimation for normal mixture distributions. *Ann. Statist.* **13** 795–800. MR790575

Hathaway, R. J. (1986). A constrained EM algorithm for univariate normal mixtures. *J. Stat. Comput. Simul.* **23** 211–230.

Hennig, C. (2003). Robustness of ML estimators of location–scale mixtures. Available at www.math.uni-hamburg.de/home/hennig/papers/hennigcottbus.pdf.

Huber, P. J. (1964). Robust estimation of a location parameter. *Ann. Math. Statist.* **35** 73–101. MR161415

Huber, P. J. (1981). *Robust Statistics*. Wiley, New York. MR606374

Keribin, C. (2000). Consistent estimation of the order of mixture models. *Sankhyā Ser. A* **62** 49–66. MR1769735

Kharin, Y. (1996). *Robustness in Statistical Pattern Recognition*. Kluwer, Dordrecht.

Lindsay, B. G. (1995). *Mixture Models*: *Theory, Geometry and Applications*. IMS, Hayward, CA.

Markatou, M. (2000). Mixture models, robustness, and the weighted likelihood methodology. *Biometrics* **56** 483–486.

McLachlan, G. J. (1982). The classification and mixture maximum likelihood approaches to cluster analysis. In *Handbook of Statistics* (P. R. Krishnaiah and L. Kanal, eds.) **2** 199–208. North-Holland, Amsterdam. MR716698

McLachlan, G. J. (1987). On bootstrapping the likelihood ratio test statistic for the number of components in a normal mixture. *Appl. Statist.* **36** 318–324.

McLachlan, G. J. and Basford, K. E. (1988). *Mixture Models*: *Inference and Applications to Clustering*. Dekker, New York. MR926484

McLachlan, G. J. and Peel, D. (2000). *Finite Mixture Models*. Wiley, New York. MR1789474

Peel, D. and McLachlan, G. J. (2000). Robust mixture modeling using the $t$ distribution. *Stat. Comput.* **10** 339–348.

Redner, R. A. and Walker, H. F. (1984). Mixture densities, maximum likelihood and the EM algorithm. *SIAM Rev.* **26** 195–239. MR738930

Rocke, D. M. and Woodruff, D. L. (2000). A synthesis of outlier detection and cluster identification. Unpublished manuscript.

Roeder, K. and Wasserman, L. (1997). Practical Bayesian density estimation using mixtures of normals. *J. Amer. Statist. Assoc.* **92** 894–902. MR1482121

Schwarz, G. (1978). Estimating the dimension of a model. *Ann. Statist.* **6** 461–464. MR468014

Tyler, D. E. (1994). Finite sample breakdown points of projection based multivariate location and scatter statistics. *Ann. Statist.* **22** 1024–1044. MR1292555

Wang, H. H. and Zhang, H. (2002). Model-based clustering for cross-sectional time series data. *J. Agric. Biol. Environ. Statist.* **7** 107–127.

Wolfe, J. H. (1967). NORMIX: Computational methods for estimating the parameters of multivariate normal mixtures of distributions. Research Memo SRM 68-2, U.S. Naval Personnel Research Activity, San Diego.

Zhang, J. and Li, G. (1998). Breakdown properties of location M-estimators. *Ann. Statist.* **26** 1170–1189. MR1635381

Zuo, Y. (2001). Some quantitative relationships between two types of finite sample breakdown point. *Statist. Probab. Lett.* **51** 369–375. MR1820795

Faculty of Mathematics-SPST
Universität Hamburg
D-20146 Hamburg
Germany
e-mail: hennig@math.uni-hamburg.de